\documentclass[11pt,letterpaper]{article}
\usepackage[margin=1in]{geometry}
\usepackage{graphicx} 
\usepackage[utf8]{inputenc} 
\usepackage{biblatex}
\usepackage{amsmath}
\usepackage{mathtools}
\usepackage{amssymb}
\usepackage{amsthm}
\usepackage{amscd}
\usepackage{epic, eepic}
\usepackage{url}
\usepackage{color}
\usepackage{makebox}
\usepackage{comment}
\usepackage{enumerate}   
\usepackage{hyperref}
\usepackage{verbatim}
\usepackage[]{todonotes}

\title{The Merino-Welsh conjecture for sparse and dense graphs}
\author{Benedek Váli\thanks{University of Szeged, Szeged, Hungary. The author is supported by the Móricz Doctoral Scholarship.}}
\date{}

\theoremstyle{plain}
\newtheorem{theorem}{\bf Theorem}
\newtheorem{lemma}[theorem]{\bf Lemma}

\newcommand\restr[2]{{
  \left.\kern-\nulldelimiterspace 
  #1 
  \littletaller 
  \right|_{#2} 
  }}

\newcommand{\littletaller}{\mathchoice{\vphantom{\big|}}{}{}{}}

\begin{document}

\maketitle
\vspace{-0.5em}\begin{abstract}
\setlength{\emergencystretch}{2em}
    The Merino-Welsh conjecture states that the Tutte polynomial $T_G(x,y)$ of a loopless and bridgeless graph $G$ satisfies the inequality $\max\{T_G(2,0),T_G(0,2)\}\geq T_G(1,1).$
    It is known to hold for several graph classes, among them sufficiently sparse and sufficiently dense graphs. 
    
    In this paper, we improve the previously known sufficient density conditions and show that the multiplicative version of the conjecture also holds for sufficiently sparse and sufficiently dense graphs. In particular, we show that for a connected graph $G$ with $m$ edges and $n$ vertices, the Merino-Welsh conjecture holds if either $m\leq1.64(n-1)$ and $G$ is loopless or $m\geq2.55(n-1)$ and $G$ is bridgeless, and the multiplicative conjecture holds if $G$ is loopless and bridgeless with $m\leq1.21(n-1)$ or $m\geq5.62(n-1)$.
\end{abstract}

\section{Introduction}
The Tutte polynomial is a two variable graph polynomial, introduced by William Tutte in 1954 \cite{tutte1954contribution}. It can be calculated recursively by $T_G(x,y)\equiv1$ when $E(G)=\emptyset$ and the following fundamental lemma.

\begin{lemma}\label{fundamental} (\cite{tutte1954contribution})
Let $G$ be a graph and $e\in E(G).$ Then:
\begin{enumerate}[(a)]
\item $T_G(x,y)=T_{G-e}(x,y)+T_{G/e}(x,y)$ if $e$ is neither a loop nor a bridge,
\item $T_G(x,y)=xT_{G-e}(x,y)$ if $e$ is a bridge,
\item $T_G(x,y)=yT_{G-e}(x,y)$ if $e$ is a loop.
\item Furthermore, if $G=G_1\cup G_2,$ where $E(G_1)\cap E(G_2)=\emptyset$ and $|V(G_1)\cap V(G_2)|\leq1,$ then $T_G(x,y)=T_{G_1}(x,y)T_{G_2}(x,y).$
\end{enumerate}
\end{lemma}

In 1999, Merino and Welsh \cite{MerinoWelsh1999} conjectured that for every loopless and bridgeless graph $G$, the Tutte polynomial $T_G(x,y)$ satisfies the inequality 
$$\max(T_G(2,0),T_G(0,2))\geq T_G(1,1),$$
and even stronger, that $T_G(x,y)$ is convex on the line segment between $(2,0)$ and $(0,2).$ Later, Merino et al. \cite{MERINO2009603} studied matroid analogues and Conde and Merino \cite{CondeMerino2009} posed additive and multiplicative variants of the conjecture, that for any bridgeless and loopless graph $G,$
$$T_G(2,0)+T_G(0,2)\geq 2T_G(1,1),$$
$$T_G(2,0)T_G(0,2)\geq T_G(1,1)^2.$$

Both the multiplicative and convex versions imply the additive version, which in turn implies the original conjecture. The values $T_G(2,0),$ $T_G(0,2)$ and $T_G(1,1)$ equal the numbers of acyclic orientations, totally cyclic orientations and maximal spanning forests of $G,$ respectively.

The original conjecture is known to hold for loopless and bridgeless threshold, complete bipartite \cite{CondeMerino2009}, subcubic, planar triangulation \cite{thomassen2010spanning} and outerplanar graphs \cite{Hansen2013}, 
and the multiplicative conjecture holds for series parallel \cite{NOBLE201424}, pathwidth $\leq3$ graphs \cite{Ok2016}. The convex version holds for minimally 2-edge-connected graphs \cite{Ok2016}.

For matroids, positive results for some matroid classes, as well as weakenings of the multiplicative conjecture were published in \cite{MERINO2009603}\cite{jackson2010inequality}\cite{CHAVEZLOMELI2011422}\cite{KNAUER201823}\cite{ferroni2023merino}\cite{BEKE2024104003}, before \cite{BEKE2024109674} showed infinitely many counterexamples to the ordinary conjecture, and hence to all versions of the conjecture for matroids. Since then, some further positive results on the conjecture and weakenings appeared in \cite{KUNG2025102866}\cite{CSIKVARI2026104402}\cite{knauer2026}, of which the first proved a sufficient density bound for matroids, namely, that if the size of a coloopless matroid of rank $r\geq4$ is $\geq r(\log_2 r+\log_2\log_2 r + \log_2\log_2\log_2 r),$ the additive conjecture holds.

Other than the above graph and matroid classes, sufficient density conditions for the ordinary conjecture were proved in \cite{thomassen2010spanning}, with improved bounds in \cite{Lin2013} and \cite{Ok2016}. Namely, the ordinary conjecture holds for connected loopless and bridgeless graphs if $m\leq 1.29(n-1),$ $m\geq 4(n-1),$ $m\geq3.58(n-1)$ with 3-edge-connectedness or $m\geq3.51n$ and the graph is simple and 3-connected with minimum degree at least 4.
The proofs share that for the sparse cases, they all show $T_G(1,1)\leq T_G(2,0)$ and for the dense cases $T_G(1,1)\leq T_G(0,2).$ More precisely, they prove upper and lower bounds that are exponential in $m$ and $n$ for the two sides of the inequality, providing a sufficient density condition for the conjecture to hold. This motivates the proof idea of obtaining an exponential bound directly for the ratios $\frac{T_G(1,1)}{T_G(2,0)}$ and $\frac{T_G(1,1)}{T_G(0,2)}$ which we achieve by deletion-contraction arguments for parallel and series classes, as in \cite{jackson2010inequality}. Our main results are the following.

\begin{theorem}\label{thm} Let $a>0,~b>1$ so that $1+ka\leq2b^k$ for all integers $k\geq1$ and $G$ be a graph with $n$ vertices, $m$ edges and $\kappa$ connected components. Then
    $$
    \begin{aligned}
    \frac{T_G(1,1)}{T_G(2,0)}\leq \frac{b^m}{a^{n-\kappa}}\qquad&\text{if $G$ is loopless, and}\\
    \frac{T_G(1,1)}{T_G(0,2)}\leq\frac{b^m}{a^{m-n+\kappa}}\qquad&\text{if $G$ is bridgeless.}
    \end{aligned}
    $$
\end{theorem}

\begin{theorem}\label{cor}
Let $G$ be a connected graph with $m$ edges and $n$ vertices. If $G$ is
\begin{enumerate}[(a)]
\item loopless and $m\leq1.64(n-1),$ then $T_G(1,1)\leq T_G(2,0).$
\item bridgeless and $m\geq2.55(n-1),$ then $T_G(1,1)\leq T_G(0,2).$
\item loopless and bridgeless with either $m\leq1.21(n-1)$ or $m\geq 5.62(n-1),$ then $${T_G(1,1)}^2\leq T_G(2,0)T_G(0,2).$$
\end{enumerate}
\end{theorem}

\section{Parallel and series classes}

All graphs considered in this paper are finite, may be disconnected and may contain both loops and multiple edges unless otherwise stated. A \textit{bridge} is an edge whose deletion increases the number of connected components of a graph, counting isolated vertices. Contraction of a set of edges deletes the edges from the graph and identifies the vertices of each connected component of the subgraph induced by the edge set while retaining loops and multiple edges among the remaining edges.
We define the parallel and series classes as in \cite{jackson2010inequality}. Call two non-loop edges $e,f$ parallel if their endpoints are the same. Then parallelism is an equivalence relation on the non-loop edges of $G$ and we denote the equivalence class of $e$ as $P(e),$ the \textit{parallel class of} $e.$ The \textit{series class} $S(e)$ is defined analogously by saying that non-bridge edges $e,f$ are equivalent if $e=f$ or $\{e,f\}$ is a cocircuit of $G$, i.e. if $\{e,f\}$ forms an edge-cut in $G.$ We say that $P(e)$ is a \textit{multiple bridge} if $e$ is a bridge in $G-(P(e)\setminus\{e\})$ and $S(e)$ is a subdivided loop if $e$ is a loop in $G/(S(e)\setminus\{e\}).$ That is, if $G-P(e)$ has more components than $G$ and if $S(e)$ is the edge-set of a circuit in $G$.

The following lemma extends and can be easily deduced from Lemma \ref{fundamental}. It is Lemma 5 from \cite{jackson2010inequality} completed with the special cases of $P(e)$ being a multiple bridge and $S(e)$ being a subdivided loop.

\begin{lemma} Let $G$ be a graph, $e\in E(G)$ and $P(e)$ and $S(e)$ be the parallel and series class of $e$ with $|P(e)|=p$ and $|S(e)|=s.$ Then
\begin{enumerate}[(a)]
\item if $P(e)$ is a multiple bridge, then
$$T_G(x,y)=xT_{G-P(e)}(x,y)+(y^{p-1}+y^{p-2}+\ldots+y)T_{G/P(e)}(x,y).$$
\item  if $P(e)$ is not a multiple bridge, then
$$T_G(x,y)=T_{G-P(e)}(x,y)+(y^{p-1}+y^{p-2}+\ldots+1)T_{G/P(e)}(x,y).$$
\item  if $S(e)$ is a subdivided loop, then
$$T_G(x,y)=(x^{s-1}+x^{s-2}+\ldots+x)T_{G-S(e)}(x,y)+yT_{G/S(e)}(x,y).$$
\item  if $S(e)$ is not a subdivided loop, then
$$T_G(x,y)=(x^{s-1}+x^{s-2}+\ldots+1)T_{G-S(e)}(x,y)+T_{G/S(e)}(x,y).$$
\end{enumerate}
\end{lemma}

Note that if $P(e)$ is a multiple bridge $T_{G-P(e)}\equiv T_{G/P(e)},$ and if $S(e)$ is a subdivided loop, $T_{G-S(e)}\equiv T_{G/S(e)}$ by e.g. Lemma \ref{fundamental} (d). Hence, for the special points (2,0), (0,2) and (1,1) we have the following.

\begin{lemma}\label{lemma} Let $G$ be a graph, $e\in E(G),$ and $P(e)$ and $S(e)$ be the parallel and series class of $e$ with $|P(e)|=p$ and $|S(e)|=s.$ Then:
\begin{enumerate}[(a)]
\item if $P(e)$ is a multiple bridge, then
$$T_G(2,0)=2T_{G/P(e)}(2,0)\quad\text{and}\quad T_G(1,1)=pT_{G/P(e)}(1,1),$$

\item  if $P(e)$ is not a multiple bridge, then
$$T_G(2,0)=T_{G-P(e)}(2,0)+T_{G/P(e)}(2,0)\quad\text{and}$$
$$T_G(1,1)=T_{G-P(e)}(1,1)+pT_{G/P(e)}(1,1)\qquad$$

\item  if $S(e)$ is a subdivided loop, then
$$T_G(0,2)=2T_{G-S(e)}(0,2)\qquad\text{and}\qquad T_G(1,1)=sT_{G-S(e)}(1,1),$$

\item  if $S(e)$ is not a subdivided loop, then
$$T_G(0,2)=T_{G-S(e)}(0,2)+T_{G/S(e)}(0,2)\quad\text{and}$$
$$T_G(1,1)=sT_{G-S(e)}(1,1)+T_{G/S(e)}(1,1).\qquad$$
\end{enumerate}
\end{lemma}

We note that if $G$ is loopless and $P$ is a full parallel class, $G-P$ and $G/P$ are also loopless by the definition of parallelism. Similarly if $G$ is bridgeless, so are $G-S$ and $G/S$ by the definition of series class. 

\section{Proof of Theorem \ref{thm}}

\begin{proof}
We prove both inequalities by induction on the number of edges. For $m=0$, both ratios equal 1 as $T_G\equiv1$, hence both inequalities hold with equality. Otherwise let $P$ and $S$ be the parallel and series class of an edge of $G$ with $|P|=p,~|S|=s.$

We first prove the first inequality. If $P$ is a multiple bridge, then
$$
\frac{T_G(1,1)}{T_G(2,0)}=
\frac{pT_{G/P}(1,1)}{2T_{G/P}(2,0)}
\leq\frac{p}2\cdot\frac{b^{m-p}}{a^{(n-1)-\kappa}}
=\frac{pa}{2b^p}\cdot\frac{b^{m}}{a^{n-\kappa}}
<\frac{1+pa}{2b^p}\cdot\frac{b^{m}}{a^{n-\kappa}}
\leq\frac{b^{m}}{a^{n-\kappa}}
$$
by Lemma \ref{lemma}, the induction hypothesis since $G/P$ is loopless and $1+ka\leq2b^k.$ 

If $P$ is not a multiple bridge, then
$$
\begin{aligned}
\frac{T_G(1,1)}{T_G(2,0)}&=
\frac{T_{G-P}(1,1)+p\cdot T_{G/P}(1,1)}{T_{G-P}(2,0)+T_{G/P}(2,0)}\\
&=\frac{\frac{T_{G-P}(1,1)}{T_{G-P}(2,0)}T_{G-P}(2,0)+p\cdot\frac{T_{G/P}(1,1)}{T_{G/P}(2,0)}T_{G/P}(2,0)}{T_{G-P}(2,0)+T_{G/P}(2,0)}\\
&\leq \frac{\frac{b^{m-p}}{a^{n-\kappa}}T_{G-P}(2,0)+p\cdot\frac{b^{m-p}}{a^{n-1-\kappa}}T_{G/P}(2,0)}{T_{G-P}(2,0)+T_{G/P}(2,0)}\\
&=\frac{b^{m-p}}{a^{n-\kappa}}\cdot\frac{1\cdot T_{G-P}(2,0)+pa\cdot T_{G/P}(2,0)}{T_{G-P}(2,0)+T_{G/P}(2,0)}
\end{aligned}
$$
by Lemma \ref{lemma} and the induction hypothesis since $G-P$ and $G/P$ are loopless. Since $G/P$ can be obtained from $G-P$ by identifying the two endpoints of $P$ and the value of the Tutte polynomial at (2,0) is the number of acyclic orientations of the graph, the weighted average of $1$ and $pa$ above is towards 1 from their arithmetic average $\frac{1+pa}2\leq b^p.$ Thus
$$
\begin{aligned}
\frac{T_G(1,1)}{T_G(2,0)}\leq \frac{b^{m-p}}{a^{n-\kappa}}\cdot\max\left\{1,b^p\right\}\leq\frac{b^m}{a^{n-\kappa}}.
\end{aligned}
$$

For the second inequality, first suppose that $S$ is a subdivided loop. Then
$$
\frac{T_G(1,1)}{T_G(0,2)}=
\frac{sT_{G-S}(1,1)}{2T_{G-S}(0,2)}
\leq\frac{s}2\frac{b^{m-s}}{a^{(m-s)-n+(\kappa+s-1)}}
=\frac{sa}{2b^s}\frac{b^{m}}{a^{m-n+\kappa}}
<\frac{b^{m}}{a^{m-n+\kappa}}
$$
by Lemma \ref{lemma}, the induction hypothesis since $G-S$ is bridgeless and $1+ka\leq2b^k.$

If $S$ is not a subdivided loop, then
$$
\begin{aligned}
\frac{T_G(1,1)}{T_G(0,2)}&=
\frac{s\cdot T_{G-S}(1,1)+T_{G/S}(1,1)}{T_{G-S}(0,2)+T_{G/S}(0,2)}\\
&=\frac{s\cdot\frac{T_{G-S}(1,1)}{T_{G-S}(0,2)}T_{G-S}(0,2)+\frac{T_{G/S}(1,1)}{T_{G/S}(0,2)}T_{G/S}(0,2)}{T_{G-S}(0,2)+T_{G/S}(0,2)}\\
&\leq \frac{s\cdot\frac{b^{m-s}}{a^{(m-s)-n+(\kappa+s-1)}}T_{G-S}(0,2)+\frac{b^{m-s}}{a^{(m-s)-(n-s)+\kappa}}T_{G/S}(0,2)}{T_{G-S}(0,2)+T_{G/S}(0,2)}\\
&=\frac{b^{m-s}}{a^{m-n+\kappa}}\cdot\frac{sa\cdot T_{G-S}(0,2)+1\cdot T_{G/S}(0,2)}{T_{G-S}(0,2)+T_{G/S}(0,2)}.
\end{aligned}
$$
by Lemma \ref{lemma} and the induction hypothesis since $G-S$ and $G/S$ are bridgeless. Since $G/S$ can be obtained from $G-S$ by identifying the endpoint pairs of the edges of $S$ and the value of the Tutte polynomial at (0,2) is the number of totally cyclic orientations of the graph, the weighted average of $sa$ and $1$ above is towards 1 from their arithmetic average $\frac{sa+1}2\leq b^s.$ Thus
$$
\begin{aligned}
\frac{T_G(1,1)}{T_G(0,2)}\leq \frac{b^{m-s}}{a^{m-n+\kappa}}\cdot\max\left\{1,b^s\right\}\leq\frac{b^m}{a^{m-n+\kappa}}.
\end{aligned}
$$
This concludes both inequalities.
\end{proof}

\section{Proof of Theorem \ref{cor}}
\begin{proof}
Let $a,b$ and $a',b'$ satisfy the condition in Theorem \ref{thm}. Then:
\begin{itemize}
\item if $G$ is loopless and $m\leq \log_{b}a\cdot(n-\kappa)$ then $T_G(1,1)\leq T_G(2,0),$
\item if $G$ is bridgeless, $a>b$ and $m\geq \frac1{1-\log_{a}b}\cdot(n-\kappa)$ then $T_G(1,1)\leq T_G(0,2),$
\item if $G$ is loopless and bridgeless, $bb'>a'$ and $m\leq \log_{\frac{bb'}{a'}}\frac{a}{a'}\cdot(n-\kappa)$ then $T_G(1,1)^2\leq T_G(2,0)T_G(0,2),$
\item if $G$ is loopless and bridgeless, $bb'<a'$ and $m\geq  \log_{\frac{bb'}{a'}}\frac{a}{a'}\cdot(n-\kappa)$ then $T_G(1,1)^2\leq T_G(2,0)T_G(0,2),$
\end{itemize}

where the last two inequalities come from
$$\frac{{T_G(1,1)}^2}{T_G(2,0)T_G(0,2)}=\frac{{T_G(1,1)}}{T_G(2,0)}\frac{{T_G(1,1)}}{T_G(0,2)}\leq\frac{b^{m}}{a^{n-\kappa}}\frac{{b'}^{m}}{{a'}^{m-n+\kappa}}=\frac{(bb'/a')^{m}}{(a/a')^{n-\kappa}}.$$

Note that the condition on $a,b$ gives $\log_ba\leq\min\{\log_b(2b-1),\log_b(b^2-\frac12)\}.$ For $b>1$, these functions are strictly decreasing in value from 2 to 1 and strictly increasing from $-\infty$ to 2, hence their minimum attains its maximum value on $b>1$ where they are equal, i.e. $b=1+\frac1{\sqrt2}.$ Since $a=1+\sqrt2,~b=1+\frac1{\sqrt2}$ satisfy $1+ka\leq2b^k~\forall k\geq1$, we get that the best constant multipliers for the first two cases that the above argument may provide are $\displaystyle\log_{1+\frac1{\sqrt2}}(1+\sqrt2)=1.648...$ and $\frac1{1-\log_{a}b}=2.543...$.

A similar argument using the constraints for $k=1,4,5$ gives the optimal constants for the multiplicative cases. Let $a,b>0$ be so that $1+ka\leq2b^k$ holds with equality for $k=4,5$ for them. Let $a'\in(4.5,5)$ be the unique solution of $a{a'}^{a'}=\left(\frac{b(a'+1)}2\right)^{a'+1}$ and $b'=\frac{1+a'}2.$ Then $a,b$ and $a',b'$ satisfy the required conditions and $(a,b,a',b')$ minimizes $\log_{\frac{bb'}{a'}}\frac{a}{a'}$ over the admissible parameters with $bb'<a'$, and $(a',b',a,b)$ maximizes $\log_{\frac{bb'}{a'}}\frac{a}{a'}$ over the admissible parameters with $bb'>a'$. 
Numerically, $(a,b,a',b')=(0.742740...,1.187042...,4.611144...,2.805572...),$ for which $bb'<a'$ and
$$\log_{\frac{bb'}{a'}}\frac{a}{a'}=a'+1=5.611...,$$
and $(a',b',a,b)=(4.611144...,2.805572...,0.742740...,1.187042...),$ for which $b'b>a$ and
$$\log_{\frac{b'b}{a}}\frac{a'}{a}=\frac{a'+1}{a'}=1.216... .$$
\end{proof}

\section*{Acknowledgements}
The proofs are due to ChatGPT 6 Astra, based on the author's motivation described above. The proofs have been verified and written up manually by the author, who takes full responsibility for the correctness of the mathematical content. ChatGPT was also used for auditing the completed manuscript.

\printbibliography
\end{document}